\documentstyle[12pt]{article}
\def\circS{{\mathop{S}\limits^\circ}}
\def\circM{{\mathop{M}\limits^\circ}}
\newcommand{\pa}{\partial}
\newcommand{\be}{\begin{equation}}
      \newcommand{\ee}{\end{equation}}
      \newcommand{\ba}{\begin{eqnarray}}
       \newcommand{\ea}{\end{eqnarray}}
\newcommand{\ban}{\begin{eqnarray*}}

\newcommand{\nin}{\noindent}
\newcommand{\f}{\mbox{\bf f}}
\newcommand{\wt}{\mbox{\bf a}}
 \newcommand{\qed}{\hspace*{\fill}$\Box$} %\rule{3mm}{3mm}\quad}

\newcommand{\sect}[1]{\section{#1} \setcounter{equation}{0}}

\begin{document}
\begin{center}
\bf{\large  On the geometry and topology of 
manifolds of positive bi-Ricci curvature
}
\end{center}
\begin{center}
Ying Shen \\Department of Mathematics, Dartmouth College\\
\end{center}
\begin{center}
Rugang Ye\\ Department of Mathematics, University of California, Santa
Barbara\\    

\&  \\

Mathematics Institute, Bochum University
\end{center}
\begin{center}
August, 1997
\end{center}

\noindent {Table of Contents} \\
\noindent  1 Introduction \\
\noindent 2 Conformal Ricci curvature \\
\noindent 3 Minimal hypersurfaces \\
\noindent 4 Topological implications \\
\noindent 5 Constructions \\
\noindent 6 Miminal surfaces of higher codimensions \\
\sect{Introduction}

In [ShY1], we introduced the concept of bi-Ricci curvature and 
initiated the study of manifolds of positive bi-Ricci curvature. 
We 
obtained a size estimate for 
stable minimal hypersurfaces in 3, 4 and 5 dimensional 
manifolds of positive bi-Ricci curvature. As a consequence, 
a homology radius estimate for the manifolds 
was derived.

In the present paper, we obtain 
a number of results regarding 
bi-Ricci curvature. 
First, we extend the size estimate in 
[ShY1] to six or less dimensional manifolds 
of positive harmonic bi-Ricci curvature.
Here,   ``harmonic bi-Ricci curvature" is 
a variant of bi-Ricci curvature. Indeed, it is 
one of the ``weighted bi-Ricci curvatures".

The above extension is achieved by 
choosing a certain conformal factor in 
our computations to be more general than 
in [ShY1]. On the other hand, we adopt a 
strategy for deriving our size estimates 
which somewhat differs from that in [ShY1].
Indeed,  
we first introduce the concept of 
``conformal Ricci curvature", and present a 
generalized Bonnet-Myers theorem about 
manifolds of positive conformal 
Ricci curvature. The desired size estimates 
for stable minimal hypersurfaces 
then follow as  
a consequence of this theorem (or its 
underlying estimate).
This way, a general geometric 
background for the topic of 
bi-Ricci curvature can be seen. 
We expect further applications of the generalized 
Bonnet-Myers theorem. (See [ShY2] for 
some closely related results on 
Riemannian and Lorentzian warped products 
and their applications to general relativity.)
On the other hand, 
we view conformal Ricci curvature 
as an interesting independent subject.

As the second topic of this paper,
we derive more geometric and 
topological 
consequences of positive (harmonic) bi-Ricci curvatures.  
On the other hand, we show that connected sums of 
manifolds of positive bi-Ricci curvatures 
admit metrics of positive bi-Ricci 
curvatures. (Note that similar 
connected sum  theorems hold for 
scalar curvature [ScY], [GL] and 
isotropy curvature [H], [MW].)
It remains to be investigated how 
a general topological classification scheme for 
manifolds of positive bi-Ricci 
curvature should proceed. But our topological 
results already provide substantial information.

Note that besides its own interest, the subject of 
positive bi-Ricci curvature has the potential 
use of providing new information 
on positive sectional curvature.  This 
is suggested e.g. by Observation 2 in the sequel.  
We hope to be able to report on this 
in the near future. 

The last topic of this paper is extension of 
our techniques to minimal surfaces of higher 
codimensions. 
Diameter estimates are derived for 2-dimensional stable minimal 
surfaces in spaces of dimensions up to 9,
and 3-dimensional stable minimal surfaces in spaces of 
dimension 5, under  
suitable positive curvature assumptions and an 
"almost flatness" condition on the 
normal bundle.  
Currently, we are trying to 
find ways to handle this almost flatness condition.  
We believe that 
further understanding of 
positive curvatures can be achieved this 
way. \\ 

\noindent {\it Acknowledgement:}  We would like to
thank 
Prof. R. Schoen for encouragement and stimulating 
discussions.

\sect{Conformal Ricci curvature }

We start with the following definition. \\

\noindent {\bf Definition 1} Let $(N, h)$ be a Riemannian 
manifold of dimension $n$ with metric $h$. For a 
positive smooth function $f$ on $N$ 
and a positive constant $\sigma$ we define 
the {\it conformal Ricci tensor} of 
$N$ associated with the conformal factor 
$f$ and the weight $\sigma$ to be 
\ba \label{2.1}
Ric^{(f, \sigma)} \equiv Ric
- \sigma (f^{-1} \Delta f) h. 
\ea
The conformal 
Ricci curvature (associated with 
conformal factor $f$ and weight $\sigma$)  in the direction of 
a unit tangent vector $v$ is then defined to 
be $Ric^{(f, \sigma)}(v, v)$.
We say that 
$(N, h)$ has positive 
conformal Ricci curvature if there is 
a conformal factor $f$ and a weight $\sigma$ 
such that 
the conformal Ricci tensor $Ric^{(f, \sigma)}$  
is positive definite.  \\

A more general version of this definition will
be given in the last section. 
The terminology in this definition is motivated by 
the transformation laws for Ricci curvature and 
scalar curvature under conformal change of 
metric.  Let $\eta$ be a positive smooth function
on $N$ and $\tilde h = \eta^2 h$ the corresponding 
conformal metric.   The said transformation laws 
are 
(see e.g. [B])

\ba \label{2.2}
\tilde Ric
= Ric
- (\Delta \ln \eta) h -(n-2)(\nabla d (\ln \eta)- d \ln \eta \otimes
d \ln \eta) - (n-2) |\nabla \ln \eta|^2 h,
\ea
and
\ba \label{2.3} 
\tilde R = \eta^{-2} (R - 2(n-1) \Delta \ln  \eta 
-(n-1)(n-2) |\nabla \ln \eta|^2),
\ea
where  $\tilde Ric = Ric_{\tilde h}, \tilde R = 
R_{\tilde h}$ ($R$ denotes the scalar curvature 
of $h$), and $\Delta$ is the trace Laplacian. 
If we set $\eta = f^{2/(n-2)}$, then (as 
is well-known)
$$
\tilde R = f^{-4/(n-2)}(R - \frac{4(n-1)}{n-2} f^{-1}
\Delta f). 
$$
This is precisely the trace of the conformal 
Ricci tensor $Ric^{(f,  4(n-1)/(n-2))}$ 
with respect to the conformal metric $f^{4/(n-2)} h$. \\

Using (2.2), we can obtain another expression for 
conformal Ricci curvature.  For given 
$f$ and $\sigma$, set $\tilde h = 
f^{2 \sigma} h$. Let $v$ be a unit tangent vector 
with respect to $h$. Consider the geodesic 
$\tilde \gamma$ with respect to $\tilde h$ which has 
$\tilde v = f^{-1} v$ as the initial tangent 
vector. Let $\gamma= \gamma(s)$ 
be the reparametrization of $\tilde \gamma$ 
which has unit speed measured in $h$. It is 
easy to see that the geodesic equation for 
$\tilde \gamma$ can be 
rewritten as follows
\ba \label{2.4}
\nabla_{\frac{\partial }{\partial s}} \frac{\partial \gamma}
{\partial s} = \sigma (\nabla \ln  f)^{\perp},
\ea 
where the $\perp$ denotes the projection to the
orthogonal complement of $\frac{\partial \gamma}{\partial s}$.
Using this equation we derive 
\ba \label{}
\nabla d(\ln  f)(\frac{\partial \gamma}
{\partial s}, \frac{\partial \gamma}
{\partial s})
=\frac {\partial^2 \ln f}{\partial s^2} - (\nabla
_{\frac
{\partial}{\partial s}} \frac{\partial \gamma}
{\partial s})
 \ln f \nonumber \\
=
\frac{\partial^2 \ln f}{\partial s^2}
 -|(\nabla \ln f)^{\perp}|^2. \nonumber  
\ea
Combining this with (2.2) (setting $\eta = f^{\sigma}$ there)
then yields 
\ba \label{2.5}
\tilde Ric(v, v)= Ric(v, v) - \Delta \ln f^{\sigma}
- (n-2) \frac{\partial^2 \ln f^{\sigma}}{\partial^2 s},
\ea      
whence
\ba \label{2.6}
Ric^{(f, \sigma)}(v, v) = \tilde Ric(v, v) 
+ (n-2) \sigma \frac{\partial^2 \ln f}{\partial s^2}
- \sigma |\nabla \ln f|^2.
\ea

\noindent {\bf Theorem 1} ({\it generalized 
Bonnet-Myers theorem})  \  {\it Let 
$(N, h)$ be connected and complete. Assume that there 
are  
a conformal factor $f$ and  weight $\sigma$, along 
with positive constants $\kappa, \varepsilon$, 
such that 
\ba \label{2.7}
Ric^{(f, \sigma)}(v, v)  + \sigma |\nabla \ln 
f |^2 - ( \frac{(n-1)}{4} + a(n) \varepsilon ) \sigma^2  |v \ln 
f|^2 \geq \kappa 
\ea
for all unit tangent vectors $v$, where $a(3)=0, a(n) =1 $ 
for $n \not = 3$. 
Then $N$ is compact.
Indeed, we have 
\ba \label{2.8}
diam (N, h) \leq \sqrt{n-1 +
\frac{(n-3)^2}{4 \varepsilon} }  \frac{\pi}{\sqrt{\kappa}}.
\ea

Consequently, the universal cover  of $N$ 
is also compact, and  hence $N$
has finite fundamental group.  }
\\

\noindent {\bf Corollary 1} {\it Assume that 
for a factor 
$f$ and a weight $\sigma < \frac{4}{n-1}$ 
($\sigma \leq \frac{4}{n-1}$ if 
$n=3$), along with a positive constant $\kappa$, there holds
\ba \label{2.9}
Ric^{(f, \sigma)} (v, v) \geq \kappa. 
\ea 
Then $N$ is compact and has 
finite fundamental group. Indeed, 
we have
\ba \label{2.10}
diam (N, h)   \leq \sqrt{n-1 + \frac{(n-3)^2}{\frac{4}{\sigma}
-n+1}} \frac{\pi}{\sqrt{\kappa}}. 
\ea
if $\sigma < \frac{4}{n-1}$, and 
\ba \label{2.11}
diam (N, h) \leq \sqrt{n-1} \frac{\pi}{\sqrt{\kappa}}
\ea
if $\sigma leq \frac{4}{n-1}$ and $n=3$.
} \\

From the viewpoint of applications, the 
following more general version of Theorem 1 
is noteworthy. 
\\

\noindent $\mbox{\bf Theorem 1}'$ {\it If we replace the constant 
$\kappa$ in the condition (2.7)
of Theorem 1 by the following quantity
\ba \label{2.12}
\frac{\kappa}{1+ dist(\cdot, p_0)^{\delta}} 
\ea
for positive constants $\kappa, \delta$ with 
$\delta <2$ and a point $p_0 \in N$, then 
we can still conclude that $N$ is compact and 
has finite fundamental group. The diameter 
of $(N, h)$ can be estimated 
as follows
\ba \label{}
diam (N, h) \leq 2 \theta^{-1}((n-1+\frac{(n-3)^2}{4\varepsilon})
\frac{\pi^2}{\kappa}), \nonumber
\ea
where $\theta$ is the increasing 
function 
\ba 
\theta(t)= \frac{t^2}{1+t^{\delta}}. \nonumber
\ea
}

Similarly, we have \\

\noindent $\mbox{\bf Corollary 1}'$ {\it If we replace 
the curvature condition in Corollary 1
by the following 
\ba \label{2.13}
Ric^{(f, \sigma)}(v, v) \geq \frac{\kappa}{1+
dist(\cdot, p_0)^{\delta}}
\ea
for positive constants $\kappa, \delta$ with
$\delta <2$ and a point $p_0 \in N$, then 
we can still conclude that $N$ is compact 
and has finite fundamental group. There holds 
\ba \label{2.14}
diam (N, h) \leq 2 \theta^{-1} 
((n-1+ \frac{(n-3)^2}{\frac{4}{\sigma}
-n+1}) \frac{\pi^2}{\kappa}) 
\ea
if $\sigma < \frac{4}{n-1}$, and
\ba \label{2.15}
diam (N, h) \leq 2 \theta^{-1}((n-1)\frac{\pi^2}{\kappa}).
\ea
if $\sigma leq \frac{4}{n-1}$ and $n=3$.
} \\

\noindent {\bf Remark 1} The classic 
Bonnet-Myers theorem follows from 
Theorem 1 by taking $f \equiv 1, 
\sigma =1$ and 
$\varepsilon$ arbitarily large.  \\

For convenience, we 
formulate two  lemmas from which 
Theorem 1 will follow. 
Consider a general manifold 
$(N, h)$, a conformal 
factor $f$, a  weight $\sigma$ and
the associated conformal 
metric 
$\tilde h = f^{2 \sigma} h$.
Let 
$\tilde \gamma$ be a unit speed geodesic in $N$ 
with repect to 
$\tilde h$ and $\gamma:
[0, 1] \rightarrow N$ its reparametrization which 
 has unit speed measured 
in $h$. 
\\

\noindent {\bf Lemma 1} {\it Assume that $\tilde \gamma$ 
minimizes length up to second order measured 
in $\tilde h$  
while its endpoints are kept fixed. 
Then we have 
\ba \label{}
(n-1)\int (\frac {d\phi}{ds})^2 + \frac {n-1}{4} 
\int \phi^2
 (\frac {d
\ln f^{\sigma}}{ds})^2 + (3-n) \int \phi
\frac {d\phi}{ds}\frac {d \ln f^{\sigma} }{ds}  \nonumber \\
\geq \int \phi^2   Ric^{(f, \sigma)}(v, v) 
+\frac{1}{\sigma} \int \phi^2 |\nabla  \ln f^{\sigma}|^2
\nonumber
\ea
for all smooth functions $\phi$ on  $[0, l]$ with
zero boundary values, where the element $ds $ in 
the integrals is omitted. } \\

\noindent {\it Proof} \  By the assumption on 
$\tilde \gamma$, there holds 

\be \label{}
(n-1)\int (\frac {d \phi}{d\tilde{s}})^2 d\tilde{s} \geq \int \phi^2 
{Ric}_{\tilde h}(\frac{\partial \tilde \gamma}{\partial \tilde s})
d\tilde{s}, \nonumber
\ee
where $\phi$ is an arbitrary smooth function on $[0, \tilde l]$
with zero boundary values. 
 
Setting $\eta = f^{\sigma}$ and applying (2.6) we 
infer

\ba 
(n-1)\int (\frac{d\phi}{ds})^2{\eta}^{-1} \geq \int \eta^{-1}
\phi^2 \{Ric^{(f, \sigma)}(\frac{\partial \gamma}{\partial s},
\frac{\partial \gamma }{\partial s})
-(n-2)\sigma \frac{\partial^2 \ln f}{\partial s^2}
-\frac{1}{\sigma} |\nabla \ln \eta|^2\}.
\nonumber 
\ea

If we replace $\phi$ by $\phi \eta^{1/2}$, 
then this inequality  becomes
\ba
(n-1)\int (\frac {d\phi}{ds})^2 + 
\frac {n-1}{4} \int \phi^2 (\frac {d
\ln \eta}{ds})^2 + (n-1)\int 
\eta^{-1}\phi \frac {d\phi}{ds}\frac {d\eta}{ds}
\nonumber \\
\geq \int \phi^2  \{ Ric^{(f, \sigma)}(
\frac{\partial \gamma }{\partial s}, \frac{\partial \gamma}{\partial s})
-(n-2) \frac {\partial ^2 \ln \eta}{\partial s^2} 
-\frac{1}{\sigma} |\nabla \ln \eta|^2\}. \nonumber
\ea
Applying the following identity we then arrive at the 
desired inequality
\ba
\int_{0}^{l}\phi^2 \frac {d^2 \ln \eta}{ds^2}=-2 \int_{0}^{l}
\phi\frac {d\phi}{ds}\frac {d \ln \eta}{ds}. \nonumber 
\ea

\qed \\

\noindent {\bf Lemma 2} {\it Let $\tilde \gamma$ satisfy the condition
in Lemma 1   
and $\gamma: [0, l] \rightarrow N$  be as before. If the 
curvature condition (2.7) holds, then  the inequality 
(2.8) holds with 
$\mbox{diam}(N, h)$ replaced by $l$.
Consequently, if the curvature condition (2.9) holds 
with $\sigma < \frac{4}{n-1}$ ($\sigma \leq 
\frac{4}{n-1}$ if $n=3$), then the inequality
(2.10) ((2.11), if $\sigma leq \frac{4}{n-1}$ and 
$n=3$) holds with $\mbox{diam}(N, h)$
replaced by $l$. } 
\\

\noindent $\mbox{\bf Lemma 2}'$  {\it Let $\tilde \gamma$
and $\gamma$ be as above. If the curvature 
condition (2.7) holds with $\kappa$ replaced by 
(2.12) (with $\delta <2$), then we have 
\ba 
l \leq \theta_{\gamma(0)}^{-1}((n-1+\frac{(n-3)^2}
{4\varepsilon})\frac{\pi^2}{\kappa}), \nonumber
\ea
where 
\ba 
\theta_{p}(t) = \frac{t^2}{1+(dist(p_0, p)+t)^{\delta}}.
\nonumber
\ea
Similarly, if the curvature 
condition (2.13) holds, then the inequality 
(2.14) (if $\sigma < \frac{4}{n-1}$) or (2.15) (if
$\sigma leq \frac{4}{n-1}$ and $n=3$)
holds with $\mbox{diam}(N, h)$ replaced 
by $l$, $\theta$ replaced by $\theta_{\gamma(0)}$,
and the factor 2 removed. } \\

With 
slight modifications, the following proof
also applies to $\mbox{Lemma 2}'$. We 
leave the details to the reader. \\

\noindent {\it Proof of Lemma 2 } \ 
Applying Lemma 1  and the assumption on 
the conformal Ricci curvature we deduce 
\ba
(n-1)\int (\frac{d \phi}{ds})^2 + (3-n) 
\int \phi \frac{d \phi}{ds} \frac{d \ln f^{\sigma}}{ds} 
- a(n) \varepsilon \int \phi^2 (\frac{d \ln f^{\sigma}}{ds})^2  \geq 
\kappa \int \phi^2. \nonumber 
\ea
If $n=3$, this implies
\ba 
(n-1)\int (\frac{d\phi}{ds})^2 \geq \kappa \int \phi^2.
\nonumber
\ea
Choosing $\phi (s) =  \mbox{sin } \frac{\pi}{l} s$ 
we then deduce 
\ba
l \leq \sqrt{\frac{n-1}{\kappa}} \pi.
\nonumber
\ea
If $n \not = 3$, we apply Yang's inequality 
and deduce
\ba
(n-1 + \frac{(n-3)^2 }{4 \varepsilon} )
\int (\frac{d\phi}{ds})^2 \geq \kappa \int \phi^2,
\nonumber
\ea
which leads to the desired  estimate.

\qed
\\

\noindent {\it Proof of Theorem 1} \
Let $p$ be a point in  $N$. Consider an arbitary $r>0$ such that
the boundary $\partial B_r(p)$ of
the geodesic ball $B_r(p)$ with respect to
$h$ is nonempty. We can find a minimizing geodesic
$\tilde \gamma$ with respect to $\tilde h$, which
runs from $p$ to $\partial B_r(p)$. Let $\gamma:
[0, l] \rightarrow N$ be its reparametrization which has
unit speed measured in $h$. Since $r \leq l$,
and $p, r$ are arbitary, Lemma 2 implies 
the desired diameter estimate. 

\qed
\\

\noindent {\it Proof of $\mbox{Theorem 1}'$} \
We choose $p=p_0$ in the above proof (replacing 
Lemma 2 by $\mbox{ Lemma 2}'$. For 
arbitary $p, q \in N$ we then use the triangular 
inequality to estimate $dist(p, q)$.

\qed 
\\

\sect{Minimal hypersurfaces} 

In this section, we apply Theorem 1 and $\mbox{Theorem 1}'$
 (or Lemma 2
and $\mbox{Lemma 2}'$) 
to study 
stable minimal hypersurfaces 
in a Riemannian manifold. 
Let $(M, g)$ be a Riemannian manifold of dimension $m=n+1$ 
with metric $g$. Consider 
a smoothly immersed minimal hypersurface $S \subset 
M$ and a bounded domain $\Omega$ in  $S$ such that 
$\bar \Omega \subset \circS$.  
The Jacobi operator  $L$
associated with the second variation of area for 
$S$ is given by  
\be \label{3.1}
L \phi= - \triangle_{S} \phi -|A|^2 \phi -Ric(\nu)\phi,
\ee
where $\nu$ denotes a unit normal of $S$,  
$A$ the second fundamental form of $S$ and $\phi$ 
a smooth function on $S$.
Let  $\lambda =\lambda_{\Omega}$ be the first eigenvalue of $L$ on 
$\Omega$ for the Dirichlet boundary value problem,
and $f$ a corresponding nonnegative first 
eigenfunction, i.e. $Lf=\lambda f$. Then $f >0$ in $\Omega$. 
We consider 
the Riemannian manifold $N=\Omega$ with 
metric  $h=g|_S$.  Choosing the eigenfunction 
$f$ as the conformal factor and an arbitary 
weight $\sigma$, we procced to 
compute conformal Ricci curvature 
$Ric^{ (f, \sigma)}_S$, where we use 
the subscript $S$ to distinguish from 
the Ricci tensor $Ric$ of $(M, g)$.

Computing at a fixed point, 
we choose an orthonormal tangent base $e_1,...,e_{n+1}$
such that $e_{n+1} = \nu$ and 
$A$ is diagonalized in the base $e_1,...,e_n$. For 
an arbitary unit tangent vector $v$ of $S$ at this point, 
we write 
$$
v = \sum_{1 \leq i \leq n} a_i e_i.
$$
Then there holds  
$$
Ric_S(v)
= \sum a_i a_j (Ric_S)_{ij}.
$$
Applying the Gauss equation and the minimality 
of $S$, we deduce 
\ba
(Ric_{S})_{ij} =\sum_{1 \leq k \leq n} R_{ikjk} + 
\sum_{1 \leq k \leq n} (A_{ij}A_{kk}-
A_{ik} A_{jk}) \nonumber \\
= R_{ij} - R_{i n+1 j n+1} - \delta_{ij} A_{ii} A_{jj}, \nonumber
\ea
where $R_{ijkl}$ is the Riemann curvature tenor of $M$.
It follows that 
\be \label{8}
Ric_S(v)
= Ric(v)-
K(v, \nu) - \sum_{1 \leq 
i \leq n} a_i^2 A_{ii }^2,
\ee
where for linearly independent tangent 
vectors $v_1, v_2$ at the same point of $M$,
$K(v_1, v_2)$ denotes the sectional curvature  of 
$M$ determined by the plane spanned by $v_1$ and $ v_2$. 
On the other hand, by the choice of $f$, we have  
\ba  
-f^{-1}\triangle_S f=Ric(\nu) + |A|^2 + \lambda.  \nonumber 
\ea
Hence we arrive at
\\

\noindent {\bf Lemma 3} {\it The conformal Ricci 
curvature of $(\Omega, g|_S)$ in the direction of $v$ 
is given by
\ba \label{}
Ric^{(f, \sigma)}_S (v, v) = Ric(v) + \sigma Ric(\nu) -K(v, \nu)
+ \sigma |A|^2 - \sum_{1 \leq
i \leq n} a_i^2 A_{ii }^2+ \sigma \lambda. \nonumber 
\ea
} \\ 

\noindent {\bf Lemma 4} {\it The following inequality holds
\ba
\sum A_{ii}^2 \geq \frac{n}{n-1} \sum a_i^2 A_{ii}^2. \nonumber
\ea
Consequently, we have
\ba
Ric^{(f, \sigma)}_S (v, v) \geq 
Ric(v) + \sigma Ric(\nu) -K(v, \nu)
+ \sigma \lambda, \nonumber
\ea
provided that $\sigma \geq (n-1)/n$.
 }\\

\noindent {\it Proof} \  By the minimality of $S$, we have
\[ A_{11}+...+A_{nn}=0. 
\]
Applying Cauchy-Schwarz inequality we then derive
\[
A_{ii}^2 =(\sum_{j \not = i} A_{jj})^2  \leq 
(n-1)\sum_{j \not = i}  A_{jj}^2 
\]
for each $i$. Consequently,
\[
a_i^2 A_{ii}^2 \leq (n-1) \sum_{j \not = i} a_i^2 A_{jj}^2.
\]
Summing over $i$ then yields 
\[
\sum a_i^2 A_{ii}^2 \leq (n-1)\sum_i A_{ii}^2 (\sum_{j \not =i}
a_j^2).
\]
Since $\sum a_i^2 =1$, adding $(n-1) \sum a_i^2 A_{ii}^2$ 
to both sides leads to the desired inequality.

\qed \\

\noindent {\bf Definition 2}  \ Let $\sigma$ be a positive 
number. For ordered orthogonal unit tangent vectors $v_1, v_2$ 
of $M$ at the same point, 
we define the {\it $\sigma$-weighted bi-Ricci curvature} of 
$M$ in 
the direction of $v_1, v_2$ to be 
\[
B_{\sigma}Rc(v_1, v_2)=Ric(v_1)+ \sigma Ric(v_2)  - K(v_1, v_2).
\] 
If $\sigma = \frac{n-1}{n} = 
\frac{m-2}{m-1}$, then $B_{\sigma}Rc$ is 
called the {\it harmonic bi-Ricci curvature}  and denoted 
by $B_HRc$. \\ 

\noindent {\bf Definition 3} 
We set 
$\partial S= \{\mbox{\it lim }
p_k:$ $ p_k \in S,$ $ \ \{p_k\}$  {\it does not converge in}
$S,$ {\it but converges in}  $M$\}. If $\partial S = \emptyset$,
we  define  $\mbox{\it diam}(S, \partial S) $
to be $\mbox{\it diam }S$. Otherwise we define it to 
be 
$\mbox{\it sup} \{ \mbox{\it dist}(p, \partial S):
p \in S\}$
\\

The following two theorems easily follow from 
Lemma 2, $\mbox{Lemma 2}'$ and Lemma 4. (If 
$\partial S = \emptyset$, then they follow from Theorem 1, 
$\mbox{Theorem 1}'$ and Lemma 4.) \\

\noindent {\bf Theorem 2} {\it Let $\sigma$ satisfy 
$(n-1)/n \leq 
\sigma < 4/(n-1) $ ($(n-1)/n \leq \sigma \leq 4/(n-1)$ 
if $n=3, i.e. \  m=4$).  
 Assume that 
$M$ satisfies the curvature condition $B_{\sigma}Rc \geq 
\kappa$ for a constant $\kappa$.  If $\bar S$ 
is noncompact, then we require $M$ to be 
complete. 
Moreover, we assume that $S$ is connected and $\lambda(S) 
> - \frac{n}{n-1}\kappa$, where
$$
\lambda(S) =   \mbox{inf } \{\lambda_{\Omega}:
\Omega
\mbox{  is a bounded domain in } S \mbox{ with }
\bar \Omega \subset \circS\}.
$$
Then 
\ba \label{12}
\mbox{diam}(S, \partial S),
dist(Q, \partial \backslash Q)
\leq c(n, \sigma)\frac{\pi}{\sqrt{\sigma \lambda(S) + \kappa}},
\ea
where $Q$ is an arbitary subset of $\partial S$ (if 
$Q$ or $\partial S \backslash Q$ is empty, then we set
$dist(Q, \partial S \backslash Q)=0$), 
and 
$$c( n, \sigma) =  \sqrt{n-1+\frac{(n-3)^2}{\frac{4}{\sigma}
-n+1}} 
$$
if $\sigma < \frac{4}{n-1}$, while 
$c(n, \sigma) = \sqrt{n-1}$ if $\sigma leq \frac{4}{n-1}$
and $n=3$.  

Consequently, if $\kappa$ is 
positive, $S$ is stable, connected  and 2-sided (i.e. 
its normal bundle is orientable), then 
\ba \label{14}
\mbox{diam} (S, \partial S), dist(Q, \partial S \backslash Q)
\leq  c(n, \sigma)\frac{\pi}{\sqrt{\kappa}}.
\ea
} \\

\noindent {\bf Remark 2} 
Note that the case of 1-weighted 
bi-Ricci curvature in dimensions $m \leq 5$ and 
the case of harmonic bi-Ricci curvature 
in dimensions $m\leq 6$ are covered by Theorem 2. 
This remark holds for all the 
results below.\\

\noindent {\bf Theorem 3} {\it Assume that $(n-1)/n
\leq \sigma < 4/(n-1)$ ( $(n-1)/n \leq
\sigma \leq 4/(n-1) $ if $m=4$),   
and there holds  
\ba
B_{\sigma}Rc \geq 
\frac{\kappa}{1+ dist(\cdot, p_0)^{\delta}}
\nonumber 
\ea
for positive constants $\kappa, \delta$ with 
$\delta <2$ and a point $p_0 \in M$. 
Let $S 
$ be a connected, 2-sided 
stable minimal hypersurface in $M$. If 
$\bar S$ is noncompact, then we require 
$M$ to be complete. Then 
\ba
diam(S, \partial S) \leq
sup\{ \theta^{-1}_p(\frac{c(m-1) \pi^2}
{\kappa}): p \in \partial S\},
\ea
\ba
dist(B, \partial S \backslash B ) \leq 
sup\{\theta^{-1}_p(\frac{c(m-1) \pi^2}
{\kappa}): p \in B\},
\ea
where $B \subset \partial S$. 
} \\

We have the following corollaries. \\

\noindent {\bf Theorem 4} {\it Assume that $M$ 
satisfy the curvature condition in Theorem 3.
Let $S$ be a complete, 2-sided stable minimal hypersurface 
in $M$. Then $S$ is compact. Moreover, its 
universal cover is also compact, and hence it 
has finite fundamental group. } \\

\noindent {\bf Theorem 5} {\it Assume that $M$ is  
an orientable,  complete Riemannian manifold 
satisfying the curvature 
condition in Theorem 3. 
Let $\Gamma$
be a smoothly embedded  
submanifold of codimension 2 in $M$.  If it is 
homologously zero and orientable, 
then it bounds a
compact smooth area-minimizing hypersurface. } \\

\noindent {\it Proof } \  
By known results in geometric measure 
theory, there is an smoothly embedded,
orientable  
area-minimizing hypersurface  with
boundary $\Gamma$. (It is obtained by 
minimizing mass for integral currents 
with $\Gamma$ as boundary.  Since integarl 
currents are oriented, the resulting 
area-minimizing hypersurface is oriented.)
Since $M$ is orientable, it follows 
that $S$ is 2-sided. Theorem 3 then implies that
it is compact. 

\qed  \\

\sect{Topological implications}

In this section we study 
topological  
implications of 
positive  bi-Ricci curvature.   
First, for the convenience of 
the reader, we recall the following definition from 
[ShY1]: \\

\noindent {\bf Definition 4} Let $M$ 
be a Riemannian manifold and {\bf V} a closed integral 
current (i.e. integer multiplicity
rectifiable current)
of codimension $2$ in $M$  
(e.g. a finite union of closed, oriented  submanifolds of
codimension 2).  (We only require  {\bf V} to have 
locally finite mass.) Assume  that {\bf V} is homologous to zero 
(in the class of integral currents of locally finite 
mass.) Then the
{\it homology
radius} of {\bf V} is defined to be
\begin{center}
$r_H(\mbox{\bf V})$ = sup  \{$r >0$, {\bf V} is 
not homologous to zero in its $r$-neighborhood\}.
\end{center}
\noindent  The homology radius $r_H(M)$ of $M$ is defined to
be the superium of $r_H(\mbox{\bf V})$ over all {\bf V}.\\ 

We shall denote the support of an integral current 
  {\bf  V} by $supp \mbox{\bf V}$. \\

\noindent {\bf Theorem 6} { \it Let $M$ be an orientable,
complete manifold of dimension $m$ such that 
$B_{\sigma}Rc \geq \kappa$ for a positive constant 
$\kappa$ and a weight $\sigma$ with $(m-2)/(m-1)\leq
\sigma < 4/(m-2)$ ($ (m-2)/(m-1) \leq \sigma \leq 4/(m-2)$
if $m=4$).
Then 
$$
r_H(M) \leq c(m-1, \sigma) \frac{\pi}{\sqrt{\kappa}}. 
$$.

More generally, let $M$ be an orientable, 
metrically complete manifold 
of dimension $m$, possibly 
with boundary, such that the abvoe curvature 
condition holds.  Let  $K$ be a subset of $M$ containing 
$\partial M$ and {\bf V} an 
integral current of codimension 2 in $M$, such that 
$[ \hbox{\bf V}] =0$ in $H_{m-2}(M, K)$
and $\hbox{dist }(supp \hbox{\bf V}, K) 
> c(m-1, \sigma)\frac{\pi}{\sqrt{\kappa}}$. Then {\bf V} 
must bound in the 
$c(m-1, \sigma)\frac{\pi}{\sqrt{\kappa}}$-neighborhood
of  $supp \mbox{\bf V}$. 
} \\

\noindent {\it Proof} \  
As in the proof of Theorem 5 (and in [Sh1Y]),
we find a mass-minimizing integral current 
{\bf W} with {\bf V} as boundary. The interior 
$S$ of $\mbox{supp } \mbox{\bf W}$ is then a smoothly 
embedded, orientable stable minimal 
hypersurface with $\partial S = \mbox{{\it supp} {\bf V}}$.
Applying Theorem 2 we then arrive at the desired 
conclusions. 
  
\qed  \\

\noindent {\bf Definition 5} Let $M$ be a noncompact 
manifold of dimension $m$.  For $1 \leq k \leq m-1$, 
let 
$H_{k}(M)_c$ denote the subgroup of
$H_{k}(M)$ generated by  chains with
compact support. A class 
$\xi \in H_k(M)_c$ is called a 
{\it free class}, provided that  it 
can be reduced to $H_{k}(M \backslash  \Omega_j)_c$
for a sequence of increasing domains $\Omega_j \subset M$ with
compact closure
such that  $M= \cup_j \Omega_j$. \\

\noindent {\bf Theorem 7} { \it Let $M$ 
be an orientable, complete Riemannian 
manifold of dimension $m$ satisfying 
the curvature condition in Theorem 3. 
Then there is no nontrivial free 
class in $H_{m-2}(M)_c$. 
Consequently, if $M$ is an orientable 
manifold of dimension $m$ with $M 
=\circM$, such that there is  a
nontrivial 
free class in $H_{m-2}(M)_c$, then
$M$ carries 
no  complete metric satisfying the 
curvature condition in Theorem 3. }
\\

\noindent {\bf Corollary 2}
{\it If $M$ is 
diffeomorphic to
the interior of a compact manifold $\bar M$ with
boundary such that the image of 
the inclusion $H_{m-2}(\partial M)
\rightarrow H_{m-2}(\bar M)$ is nonzero, then 
$M$ carries no complete metric satisfying 
the curvature condition in Theorem 3. 
} \\

\noindent {\bf Corollary 3} 
{\it If $M= M' \times S^1$ or $M = M' \times \mbox{\bf
R}^1$ such that $M'$ is an orientable 
$\mbox{\bf Z}$-periodic
manifold, then $M$ carries no complete 
metric satisfying the curvature condition of 
Theorem 3.  } \\

\noindent {\it Proof of Theorem 7}  \  
Assume the contrary. Then  
we can find integral currents with 
compact support  $\mbox{\bf V}$ and $ \mbox{\bf V}_k, 
k=1, 2,...,$ such that 
neither of them is
homologous to zero, each $\mbox{\bf V}_k$ 
is homologous to {\bf V}, 
and  
$\mbox{\it dist}(\mbox{\it supp}
 \mbox{\bf V},
\mbox{\it supp} \mbox{\bf V}_k) \rightarrow \infty$. 
For each $k$,
we can find a mass-minimizing integral current
$\mbox{\bf W}_k$ with boundary $\mbox{\bf V} - \mbox{\bf V}_k$.
One connected 
component $S_k$ of the interior of $ supp \mbox{\bf W}_k$ 
is a smoothly embedded, one-sided 
stable minimal hypersurfaces such that 
$\partial S_k = Q_k \cup T_k$ with $Q_k \not 
= \emptyset, T_k \not = \emptyset$ and
$Q_k \subset \mbox{\it supp} \mbox{\bf V}, 
T_k \subset \mbox{\it supp} \mbox{\bf V}_k$. 
Applying (3.6) in Theorem 3 we then arrive 
at a contradiction.

\qed \\

\noindent {\bf Definition 6} We say that  $M \in \cal B$,
provided that $M$ is a compact orientable manifold 
of dimension $m$ such that 
\\
(1) if $\tilde M$ is a finite covering of $M$, then
every nontrivial class in  $H_{m-1}(M)$ can be
represented by a finite disjoint union of
embedded, compact orientable hypersurfaces (with
multiplicities) which have finite fundamental
groups,
\\
(2) if $\tilde M$ is a noncompact covering of
$M$, then there is no nontrivial free class in
$H_{m-2}(\tilde M)_c$.  \\

\noindent {\bf Theorem 8} {\it Let $M$ be a compact, 
orientable $m$-dimenional 
manifold of positive harmonic bi-Ricci 
curvature with $m \leq 6$ (or positive 
1-weighted bi-Ricci curvature with 
$m \leq 5$). Then $M$ belongs to the class $\cal B$. } \\

Since positive bi-Ricci curvatures imply positive 
scalar curvature, we note the following  additional
topological implications of positive bi-Ricci
curvatures for a compact manifold $M$ of dimension
$m$: \\
(1) if $m \leq 7$, then $M \in 
{\cal C}_m$, where ${\cal C}_m$ is the 
class of manifolds introduced in [ScY2], \\
(2) if $M$ is spin, then its $\hat A$-genus
and Hitchin invariant vanish,
\\
(3) if $m=4$ and $b_2^+(M) >1$ 
for one orientation of $M$, then
the Seiberg-Witten invariants
of $M$ (for that 
orientation) vanish (if $b_2^+(M) =1$,
then the Seiberg-Witten invariants of 
$M$ corresponding to the 
chamber of the given metric on $M$ 
vanish).  \\

The above results, 
along with the constructions of manifolds of 
positive bi-Ricci curvature in the 
next section, provide a basis for 
topological classification of 
manifolds of positive bi-Ricci curvatures 
in dimensions $m \leq 6$. 
\\

\sect{Constructions}

We first make a few simple observations about the
concept of bi-Ricci curvatures. \\

\noindent {\bf Observation 1} \  Summing over an orthonormal
base we see that positive ($\sigma$-weighted)
bi-Ricci curvature implies positive scalar curvature.
On the other hand, positive sectional curvature
obviously implies positive bi-Ricci curvature. \\

\noindent {\bf Observation 2} \  
If $M$ admits a metric
with positive bi-Ricci curvature
and positive Ricci curvature (in particular, if $M$
admits a metric with positive sectional curvature), then
$M \times S^1$ and $M \times \mbox{\bf R}^1$
admits a metric with positive bi-Ricci curvature.
Taking known examples of manifolds with
positive sectional curvature, we then obtain
manifolds which admit metrics with positive
$\sigma$-weighted bi-Ricci curvature for all
$\sigma$, but admit no metric with positive
Ricci curvature.  \\

Next we present our result about connected sums.
\\

\noindent {\bf Theorem 9} {\it Let $M_1, M_2$ be two 
manifolds of the same dimension $m \geq 3$ such that 
they  both admit metrics of positive  $\sigma$-weighted
bi-Ricci curvature for the same weight 
$\sigma > 0$. Then their connected sum $M_1  \#  M_2$ 
also admits  
metrics of positive $\sigma$-weighted 
bi-Ricci curvature. 
Moreover, if $M_1, M_2$ both admit metrics whose 
$\sigma$-weighted bi-Ricci curvature 
is positive and satisfies a large scale condition 
(such as a growth condition), then $M_1 \# M_2$ 
also admits such metrics.  }
\\

This theorem clearly follows from Lemma 6 below.
\\

\noindent {\bf Lemma 5} {\it Let $(M, g)$ be a Riemannian 
manifold  of dimension $m \geq 3$. 
 Let $p \in M$ and 
$B$ be a compact geodesic ball centered at $p$ on which 
the $\sigma$-weighted bi-Ricci curvature for a weight 
$\sigma > 0$ is  positive.
Then there is a 
positive number $\rho_0$ with the following 
property. For arbitary positive numbers $\rho 
\leq \rho_0$ and $\delta$, there is  a positive 
smooth function $\eta$ on $M \backslash \{p\}$
such that $\eta|_{M \backslash B} \equiv 1$ and 
the conformal metric $\eta^2 g$ has positive 
$\sigma$-weighted bi-Ricci curvature on $B$. Moreover, 
measured in the $C^2$ norm and 
near $r =0$ ($r$ is the distance to 
$p$), the conformal metric $\eta^2 g$ 
is uniformly within $\delta$ distance 
from the standard
product metric on $\mbox{\bf  R}^+
\times S^{m-1}(\rho)$.  
Here $S^{m-1}(\rho)$
denotes the round $(m-1)$-sphere of radius $\rho>0$.  
} \\

\noindent {\bf Lemma 6} {\it  Let $(M, g)$, $\sigma$, 
$\kappa$, $p$,
$B$ and 
$\rho_0$ be as above. For arbitary 
positive number $\rho \leq \rho_0$, 
there is a smooth metric on 
$M \backslash \{p\}$, such that 
its $\sigma$-weighted 
bi-Ricci curvature is positive on $B\backslash \{p\}$,
and near $r=0$ it coincides with the standard product metric 
on $\mbox{\bf  R}^+
\times S^{m-1}(\rho)$. }\\

\noindent {\it Proof of Lemma 5} \  We follow a scheme   
in [MiW], where a connected sum theorem about positive
isotropy 
curvature is proved.  Set $\tau= \ln \eta$. We would like 
to determine $\tau$ so that the conformal 
metric $\tilde g =  \eta^2 g$ has the desired 
properties.   We choose $\tau=\tau(r)$. (We will make 
sure that $\tau$ is smooth.) 
Then  
$$
d \tau= \tau' dr, \nabla d \tau= \tau'' dr \otimes dr + \tau' \nabla d r.
$$

Consider a unit tangent vector $v =
a e_1 + b \frac{\partial}{\partial r}$, where $e_1$ 
is a unit tangent vector perpendicular to  
$\frac{\partial}{\partial r}$.  By (2.2) and 
the above formulas, we have
\ba
\tilde Ric(v) = Ric(v) - \tau''(1 + (m-2)b^2)  - \tau' \Delta r 
- (m-2) a^2 \tau' \nabla r (e_1, e_1)
\nonumber \\
-(m-2) a^2 (\tau')^2. \nonumber 
\ea
But $\nabla d r (e_1, e_1) = 1/r + O(r), 
\Delta r = (m-1)/r + O(r)$ for $r$ small.  Hence we have 
(for $r$ small) \\
\ba
\tilde Ric(v) = Ric(v) - \tau''(1 +(m-2) b^2)
\nonumber \\
-\tau' ((m-1) 
+ (m-2)a^2) ( \frac{1}{r} + O(r))
- (m-2)  a^2 (\tau')^2.
\ea
To proceed, we choose $\tau$ such that 
\ba
\tau' = - \frac{\phi}{r}
\ea
for a function $\phi$ which is to be determined. 
Then 
$$
\tau''=- \frac{\phi'}{r} + \frac{\phi}{r^2}.
$$
It follows that 
\ba
\tilde Ric(v) = Ric(v) + \frac{\phi}{r^2} 2(m-2)
(a^2  + O(r^2)) \nonumber \\
- (m-2) a^2 \frac{\phi^2}{r^2} +\frac{\phi'}{r}
(1+(m-2)b^2).
\ea
If $w = a_1 e_2 + b_1 \frac{\partial }{\partial r}
$ is another unit tangent vector such that 
$e_2 \perp \frac{\partial} 
{\pa r}$ and $v \perp w$, then we have
\ba 
\tilde Ric(v) + \sigma \tilde Ric(w) =
Ric(v)+ \sigma Ric(w) + \frac{\phi}{r^2}2 (m-2)
( a^2 + \sigma a_1^2 + O(r^2)) \nonumber \\
- (m-2) (a^2 +\sigma a_1^2) \frac{\phi^2}{r^2} +\frac{\phi'}{r}
(1+(m-2)(b^2+ \sigma b_1^2).
\ea
On the other hand, the following law holds for 
Riemann curvature tensor (see [B]),
\ba
 e^{-2\tau} \tilde Rm = Rm - g \Diamond (\nabla d \tau 
- d\tau \otimes d\tau + \frac{1}{2}|\nabla \tau|^2 g),
\ea
where for symmetric 2-tensors $A, B$, 
\ba
(A \Diamond B)(v_1, w_1, v_2, w_2) = 
A(v_1, v_2) B(w_1, w_2) + A(w_1, w_2) B(v_1, v_2)
\nonumber \\
-A(v_1, w_2) B(w_1, v_2) -A(w_1, v_2) B(w_1, v_2).
\nonumber
\ea
Simple computations then yield for 
$\tilde v = e^{-\tau} v, \tilde w = e^{-\tau} w$
\ba
\tilde Rm(v, w, \tilde v, \tilde w) =
K(v, w) -\nabla d \tau(v, v) -\nabla d \tau(w, w)
\nonumber \\
+ |v \tau|^2 + |w \tau|^2 - |\nabla \tau|^2.
\ea
Consequently, 
\ba
\tilde Rm (v, w, \tilde v, \tilde w) =
K(v, w) + \frac{\phi}{r^2}2 (a^2 +a_1^2 -1)(1 + O(r^2))
\nonumber \\
- \frac{\phi^2}{r^2}(a^2 + a_1^2) + \frac{\phi'}{r}( b^2 + b_1^2).
\ea
Combining this with (5.4) we then deduce
\ba
e^{2 \tau} B_{\sigma}\tilde Rc (\tilde v, \tilde w) 
= B_{\sigma}Rc(v, w) + \frac{\phi}{r^2}2( 1-a_1^2+ (m-3) a^2 + (m-2)\sigma
a_1^2 
+O(r^2)) \nonumber \\
 - \frac{\phi^2}
{r^2}((m-2)( a^2 + \sigma a_1^2) - a^2 - a_1^2) 
+\frac{\phi'}{r}(1-b_1^2+ (m-3) b^2 + (m-2)\sigma b_1^2)).
\nonumber 
\ea
Since $1-a_1^2 +(m-3)a^2 + 
(m-2)\sigma a_1^2  \geq min\{1, \sigma\} (1-a_1^2 +
(m-2) a_1^2) \geq min\{1, \sigma\}$, this 
formula implies  that 
\ba
e^{2 \tau} B_{\sigma}\tilde Rc (\tilde v, \tilde w)
\geq    B_{\sigma}Rc(v, w) + \frac{1}{r^2}(c_0- c_1{\phi^2})
{r^2}+ c_1\frac{\phi'}{r}, 
\ea
provided that $r \leq r_0$ for a positive number 
$r_0$ and 
$\phi' \leq 0$, where $c_0 = min\{1, \sigma\}/2$ 
and $
c_1 = m+ (m-2) \sigma $.  
Let $\kappa$ denote the minimum of $B_{\sigma}Rc$ 
on $B$. We choose $r_0$ such that $B_{r_0} \subset 
B$.
Then we seek a nonincreasing $\phi$ such that
\ba
\kappa + \frac{1}{r^2}\phi(c_0-c_1 \phi)
+c_1\frac{\phi'}{r}> 0.
\ea
Consider the 
change of variables $r = r_0 e^{-t}$ mapping 
the interval $[0, \infty)$ onto $(0, r_0]$.
Then  (5.9) is equivalent to
\ba
c_1 \psi' < \kappa r_0^2 e^{-2t} +  
\psi (c_0- c_1 \psi)
\ea
with $\phi(s)= \psi(r_0 e^{-t})$ (we require  $\psi' \geq 0$).
Rescaling, we 
can reduce this to
\ba
\psi' < c_1^{-1} \kappa r_0^2 e^{-2t} +  
\psi (c_1^{-1}c_0 - \psi)
\ea

Observe 
that the solutions of the 
ODE 
$
\psi'=\psi(c_1^{-1} c_0- \psi)
$
are given by
$$\psi_c(t) =
\frac{ c_1^{-1}c_0 c e^{c_1^{-1}c_0 t}}{1+c e^{c_1^{-1}c_0 t}}
$$
for arbitary constants $c$. 
For each positive number $t_1  >  2 \ln  2$, choose a
smooth non-decreasing function $ \beta=\beta_{t_1}$ 
such that $\beta(t) \equiv 0 $ for $t \leq \ln 2, \beta(t) 
= e^{c_1^{-1}c_0 t}$ for $ 2 \ln 2 \leq t \leq t_1$ and 
$\beta(t) \equiv e^{c_1^{-1}c_0(t_1 +1)} $ for $t \geq t_1+1$. 
We set
$$
\tilde \psi_{c}(t) = \frac{c_1^{_1}c_0 c \beta(t)}{
1+c \beta(t)}.
$$
Since $\kappa$ is positive,
it is easy to see that for each $t_1$, if we choose 
$c= c(t_1)$ to be small enough, then $\tilde \psi$
will satisfy the inequality (5.11) (and hence 
(5.10) after scaling). 
The corrsponding 
$\phi$ then satisfies (5.9), equals $1$ near $0$,
and equals $0$  for $r \geq r_0/2$.  
For this choice of $\phi$, we have $B_{\sigma}\tilde Rc 
> 0$ on $B$. Furthermore, by (5.2), there holds 
\ba
e^{2\tau(r)} = e^{2 \tau(r_1)} \frac{r^2_1}{r^2},
\ea
\ba
\tau(r_1)= \int_{r_1}^{r_0} \frac{\phi}{r},
\ea
where $r_1 = r_0 e^{-(t_1+1)}$.  
Since 
$g$ is asymptotically euclidean 
at $p$, (5.12) implies that 
near $r =0$ or $s= \infty$ under 
the tansformation $s= -\ln r$, 
the conformal metric $\tilde g$ approaches 
the 
metric 
$e^{2 \tau(r_1)} \bar r_1^2
(ds^2 + d \omega^2)$, where 
$d \omega^2$ denotes the metric of the unit sphere.
But the latter is equivalent to the 
product metric on 
$\mbox{\bf R}^+ \times S^{m-1}(\rho)$,
where $\rho = e^{\tau(r_1)} r_1$. One easily 
works out  the explicit formula for $\rho$ and 
shows that it can be made arbitarily small 
by choosing $t_1$ large.  

\qed
\\

\noindent {\it Proof of Lemma 6} \  We interpolate 
the product metric on $\mbox{\bf R}^+ \times 
S^{m-1}(\rho)$ with $\tilde g$. 

\qed

\sect{Minimal surfaces of higher codimensions}

First, we present an extension of the concept of 
conformal Ricci curvature and the generalized Bonnet-Myers
theorem.  Consider a Riemannian manifold $(N, h)$ of 
dimension $n$.\\

\nin {\bf Definition 7}  Let $k$ be a natural number, 
$f_1,...,f_k$ positive smooth functions on $N$ and 
$a_1,...,a_k$ positive numbers. The conformal 
Ricci tensor of $(N, h)$ associated with the 
conformal factor $\mbox{\bf f}=(f_1,...,f_k)$
and weight $\mbox{\bf a} = (a_1,...,a_k)$ 
is defined to be 
$$
Ric^{(\f, \wt)} = Ric- (\sum f_i^{-1} \Delta f_i)h.
$$
\\

For a conformal factor $\f$ and a weight $\wt$, we 
consider the conformal metric $\tilde h = f_1^{2a_1}
\cdot \cdot \cdot f_k^{2a_k} h$. Using (2.5) with
$f^{\sigma}$ replaced by $\eta =f_1^{a_1} \cdot \cdot \cdot 
f_k^{a_k}$  we deduce
\ba
Ric^{(\f, \wt)}(v, v)= \tilde Ric(v, v) +(n-2)
\frac{\partial^2 \ln \eta}{\partial s^2}
-\sum a_i |\nabla \ln f_i|^2
\ea
for unit tangent vectors $v$, where 
we use a curve 
$\gamma(s)$ in a similar way to 
the context of (2.5).  
Next let $\gamma: [0, l] \rightarrow
N$ be a unit speed (measured in 
$h$) curve which
minimizes length in $\tilde h$
up to second order. 
Arguing as in the proof of Lemma 1 we 
infer 
\ba
(n-1) \int (\frac{d \phi}{ds})^2 + 
\frac{n-1}{4} \phi^2 (\sum \frac{d \ln f_i^{a_i} }{ds})^2 
+ (3-n) \int \phi \frac{d\phi}{ds} \frac{d \ln \eta}
{ds}  \nonumber \\
\geq \int \phi^2 Ric^{(\f, \wt)}(\frac{\partial \gamma}{\partial 
s}, \frac{\partial \gamma}{\partial s}) + \int \phi^2 
(\sum \frac{1}{a_i} |\nabla \ln f_i|^2).
\ea
Now, the arguments in the proof of Lemma 2 and 
Cauchy-Schwarz inequality  lead 
to \\

\nin {\bf Lemma 7} {\it 
Assume $a \equiv \sum a_i < \frac{4}{n-1}$
($ a \leq \frac{4}{n-1} $ if $n=3$) and 
$Ric^{(\f, \wt)} \geq \kappa$ for a positive 
constant $\kappa$. Then we have
\ba
l \leq \sqrt{n-1 + \frac{(n-3)^2}{\frac{4}{a}
-n+1}} \frac{\pi}{\sqrt{\kappa}}
\ea
if $a <   \frac{4}{n-1}$, and  
\ba
l \leq  \sqrt{n-1}\frac{\pi}{\kappa}
\ea
if $a \leq  \frac{4}{n-1} $ and $n=3$.
}\\

Of course,  all the statements in Lemma 2, $\mbox{Lemma 2}'$,
Theorem 2 and $\mbox{Theorem 2}'$ can be extended to the present 
situation. We leave the details to the reader.

Next let $(M, g)$ be a Riemannian manifold of dimension 
$m > n$ and $S \subset M$ an immersed stable minimal surface
of dimension $n$.  
We have the following formula for
the second variation of area [Si]
\ba  
I(X)
= \int_S(\sum|(\nabla_{e_i}X)^{\perp}
|^2 - \sum<X, A(e_i, e_j)>^2 - \sum Rm(e_i, X, e_i, X)),
\ea
where $X$ is a normal vector field and $\{e_i\}$ a local
orthonormal tangent  base of $S$.  Consider 
a bounded domain $\Omega$ such that $\bar \Omega 
\subset \circS$ and there is a 
smooth unit normal vector field $\nu$ on 
$\bar \Omega$.  Choosing $X= \phi \nu$ in (6.1) we 
obtain 
$$
\int_S ( |\nabla \phi|^2 + \phi^2 |\nabla^{\perp}
\nu |^2 - \phi^2 |A_{\nu}|^2 - \sum Rm(e_i, \nu, e_i, \nu)
\phi^2) \geq 0
$$
for all smooth functions $\phi$ 
with compact support in $\Omega$, where $\nabla^{\perp}$
denotes the normal connection and 
$A_{\nu} = A\cdot \nu$.  Let 
$f$ be a positive (in $\Omega$) 
first eigenfunction of the operator 
$L \phi = -\Delta_S \phi + |\nabla^{\perp} \nu|^2 
\phi - |A_{\nu}|^2 \phi - \sum Rm(e_i, \nu, e_i, \nu) \phi$.
Then $L f  \geq 0$.  Consider the Riemannian 
manifold $N = \Omega$ with metric $g|_S$. 
We compute the conformal Ricci curvature 
$Ric^{(f, \sigma)}(v, v)$ for a weight $\sigma$ and 
a unit tangent vector $v$.  
Choose local orthonormal normal frame 
$\nu_{\alpha}$. Then we have by the Gauss equation
$$
Ric_S(v) = \sum_{i} K(v, e_i) 
+ \sum_{ \alpha} \sum_i (A_{\nu_{\alpha}}(v, v)
 \cdot A_{\nu_{\alpha}}(e_i, e_i) -  A_{\nu_{\alpha}}(v, 
e_i) \cdot A_{\nu_{\alpha}}(v, e_i) ).
$$
For each fixed $\alpha$, 
the quantity $\sum_i (A_{\nu_{\alpha}}(v, v)
 \cdot A_{\nu_{\alpha}}(e_i, e_i) -  A_{\nu_{\alpha}}(v, 
e_i) \cdot A_{\nu_{\alpha}}(v, e_i) ) =
-\sum_i  A_{\nu_{\alpha}}(v, 
e_i) \cdot A_{\nu_{\alpha}}(v, e_i) $ (by the 
minimality)  is independent 
of the choice of $e_i$. We switch to a base $\{e_{i, \alpha}\}$
to diagonalize $A_{\nu_{\alpha}}$.  In this base, we have 
$v = \sum_i a_{i, \alpha} e_{i, \alpha}$ for suitable $a_{i, \alpha}$.
It follows that
\ba
Ric_S(v) = \sum_{i} K(v, e_i)
-\sum_{\alpha} \sum_i a_{i,\alpha}^2 A_{\nu}(e_{i, \alpha},
e_{i, \alpha})^2.
\ea
Applying Lemma 4 we arrive at
\\

\noindent {\bf Lemma 8} {\it  There holds \\
\ba
Ric^{(f, \sigma)}(v, v) \geq \sum K(v, e_i) 
\sigma \sum K(\nu, e_i) \nonumber \\
 + \sigma |A_{\nu}|^2 - \frac{n-1}{n}|A|^2 
- \sigma |\nabla^{\perp} \nu|^2.
\ea
}\\

Now we assume that there is a global smooth orthonormal frame 
$\{\nu_{\alpha}\}$ for the normal bundle over $\bar \Omega$. 
For each $\nu = \nu_{\alpha}$ we choose a positive 
first eigenfunction $f_{\alpha}$. Then we obtain \\

\nin {\bf Lemma 9} {\it For $\f = (f_1,...,f_{m-
n})$ 
there holds
\ba
Ric^{(\f, \wt)}(v, v) \geq \sum_i K(v, e_i)
+ \sum_{\alpha} \sigma_{\alpha} 
(\sum_i K(\nu_{\alpha}, e_i)) \nonumber \\
+ \sum_{\alpha} \sigma_{\alpha} |A_{\nu_{\alpha}}|^2
- \frac{n-1}{n} |A|^2 - \sum_{\alpha} 
 \sigma_{\alpha} |\nabla^{\perp} \nu_{\alpha}|^2.
\ea
} \\

We choose $a_{\alpha} = \sigma$ for 
all $\alpha$ and a positive $\sigma $. 
We can apply Lemma 7 and Lemma 9 to the situation
$(n-1)/n < \sigma < 4/((m- n)(n-1))$ (or 
$(n-1)/n < \sigma \leq 4/{(m- n)}(n-1)$ if $n=3$). There are 
two possible cases (besides the already 
treated case $m=n+1$):
{\it Case 1} \ $n=2$, $4 \leq m \leq 9$ and 
$1/2 < \sigma < 4/{(m-n)}$); \\
{\it Case 2} \ $n=3, m=5$ and 
$ 2/3 < \sigma \leq  1$. \\

\noindent {\bf Definition 8} Let $p \in M $ and 
$V$ be a nontrivial proper subspace of $T_pM$.  
For a unit tangent vector $v \in V$  and 
a positive weight $\sigma$
we set 
$$
K_{\sigma}(v, V) = \sum_i K(v, e_i) + \sigma \sum_{i, \alpha}
K(\nu_{\alpha}, e_i),
$$
where $\{e_i\}$ denotes an orthonormal base of 
$V$ and $\{\nu_{\alpha}\}$ an orthonormal base of 
$V^{\perp}$. (This expression is indepedent of 
the choice of the bases.)  
\\

We leave the version of the following result 
based  on
$\mbox{Lemma 2}'$ (instead of Lemma 2) to the reader.
\\

\noindent {\bf Proposition 1} {\it Assume that we are in one 
of the above two cases. 
Moreover, assume that 
there is a positive constant $\kappa$ 
such that $K_{\sigma}(v, V) \geq \kappa$ for 
all $V$ and $v \in V$ with $dim V = n$. 
Let $S'$ be a connected domain of $S$ 
such that the normal bundle of $S'$ is flat 
and trivial. 
More generally, assume that the normal 
bundle of $S'$ admits a global orthonormal 
frame 
$\{\nu_{\alpha}\}$ which is ``almost flat" in the sense that 
\ba
\sum |\nabla^{\perp} \nu_{\alpha}(p)|^2  \leq (\sigma -
\frac{n-1}{n}) |A|^2
+ \epsilon (1+ \mbox{ dist}^2(p, \partial S')) 
\ea
for a positive constant $\epsilon$ with 
\ba
\epsilon < \frac{\kappa^2}{2(\kappa +
1 + \frac{(m-2) \sigma}{4 -( m-2) \sigma} \pi^2)}
\ea
in Case 1, or  
\ba
\epsilon  < \frac{\kappa}{8 \pi + 2\kappa}
\ea
in Case 2.

Then we have for $Q \subset \partial S'$ 
\ba
diam(S', \partial S'), dist(Q, \partial S'
\backslash Q) \leq  \sqrt{1 + \frac{(m-2) \sigma}
{4-(m-2) \sigma}} \frac{\sqrt{2}\pi}{\sqrt{\kappa-2 \epsilon}}
\ea
in Case 1, or
\ba
diam(S', \partial S'), dist(Q, \partial S' \backslash
Q) \leq   \frac{2\pi}{\sqrt{\kappa -2 \epsilon}}
\ea
in Case 2.
} \\

\noindent {\it Proof} \  We prove the estimate for 
$diam(S', \partial S')$ in Case 2. 
The other estimate and case
can be handled in a similar way.  First assume 
that $S'$ is bounded and $\bar S' \subset 
\circS$.  Consider the positive first 
eigenfunctions $f_{\alpha}$ on $S'$ corresponding to 
$\nu_{\alpha}$ and the associated conformal metric 
$\tilde g = f_1^{2\sigma} \cdot \cdot \cdot 
f_{m-n}^{2 \sigma} g_S$. For 
small  $t>0$,  consider an arbitary connected 
subdomain $\Omega $ 
in $S'$ 
such that $\bar \Omega$ is contained in the interior of 
$ S'$ and 
$dist(p, \partial S') \leq t$ for all 
$p \in \partial \Omega$.  

Set $ F = \{ p \in \bar \Omega: dist(p, \partial \Omega)
\leq \sqrt{\frac{ \kappa}{2 \epsilon} -t}.$   We claim 
that 
$F$ is open in $\bar \Omega$. It suffices 
to show that no point in $
\bar \Omega$ has distance $\sqrt{\frac{ \kappa}{2 \epsilon}} -t$
to $\partial \Omega$. Assume the contrary and 
let $p \in \bar \Omega$ has this distance to $
\partial \Omega$. Let $\tilde \gamma$ be 
a minimizing geodesic from $p$ to $\partial 
\Omega$ with respect to the metric $\tilde g$ and 
$\gamma:[0, l] \rightarrow \bar \Omega'$
its reparametrization with unit tangent 
vector measured in $g_S$. By Lemma 7 
and the almost flatness assumption, we deduce 
($\wt = (\sigma,...,\sigma)$)
$$
Ric^{(\f, \wt)}(\frac{\partial \gamma}{\partial s},
\frac{\partial \gamma}{\partial s}) \geq \frac{\kappa}{2}
- t \epsilon. 
$$
Applying Lemma 7 we then infer 
$$
l \leq  \frac{2 \pi}{\sqrt{\kappa-2 \epsilon }}
< \sqrt{\frac{ \kappa}{2 \epsilon}} - t \epsilon,
$$
provided that $t$ is small 
enough. This is a contradiction. Since $F$ is 
also closed, it follows that $F = \Omega$. Approximating $S'$ by 
$\Omega$ we then obain the desired estimate.  If $S'$ is 
not bounded, we can approximate it by bounded 
domains.  

\qed 
\\

\noindent {\bf Remark 3} \ We can replace the term 
$dist^2(p, \partial \Omega)$ in (6.9) by
other increasing 
functions of $dist(p, \partial \Omega)$. 
Of ocurse, the requirements (6.10) and (6.11)
on $\epsilon$ 
will change accordingly.  
\\

\nin {\bf Remark 4} \ The triviality condition 
on the normal bundle in Proposition 1 is not 
so restrictive as it might appear at the first 
glance. Indeed, for example,
if $S$ is topologically a sphere,
then $S'= S \backslash \{p\}$ for any $p \in S$ 
satisfies this condition. In contrast, the 
almost flatness condition is more serious. 
We would like to point out that 
if the normal curvature is 
small, then it is easy to construct 
almost flat frames in an injectivity 
domain of $S$.  On the other hand, by
the Ricci 
equation, small normal curvature follows 
from a suitbale pinching condition on 
the curvature tensor of $M$. As a consequence, 
some local size estimates for $S$  can be derived
under 
such a pinching condition. For example, 
we can deduce an upper bound for the conjugate radius 
at every point. By the Gauss equation, this implies 
an estimate for the second fundamental form in the 
case $dim S =2$.    Details wil be given   
elsewhere.

\end{document}